\documentclass[12pt,reqno,twoside]{amsart}
\usepackage{amsthm,amsmath,amssymb,mathrsfs,amsbsy}
\usepackage{mathabx} 
\usepackage{graphicx,epsfig,wrapfig,sidecap,subfig}
\usepackage[all]{xy}
\usepackage{proof}
\usepackage{xcolor}
\usepackage{hyperref}
\usepackage{bm}
\usepackage[capitalize]{cleveref}
\usepackage{listings}
\usepackage[normalem]{ulem}

\usepackage[a4paper, hmargin=2.5cm, vmargin={3.1cm, 3.1cm}]{geometry}

\usepackage{color}

\title[Discrepancy and rectifiability of almost linearly repetitive Delone sets]{Discrepancy and rectifiability of almost linearly repetitive Delone sets}
\author{Yotam Smilansky}
\address{Yotam Smilansky\newline\indent Department of Mathematics, Rutgers University, NJ, USA. \newline\indent {\tt yotam.smilansky@rutgers.edu}}
\author{Yaar Solomon}
\address{Yaar Solomon\newline\indent Department of Mathematics, Ben-Gurion University of the Negev, Israel.\newline\indent {\tt yaars@bgu.ac.il}}

\newcommand{\N}{{\mathbb{N}}}
\newcommand{\Z}{{\mathbb{Z}}}
\newcommand{\Q}{{\mathbb {Q}}}
\newcommand{\R}{{\mathbb{R}}}

\newcommand{\X}{{\mathbb{X}}}

\newcommand{\BB}{\mathcal{B}}
\newcommand{\CC}{\mathcal{C}}

\newcommand{\TT}{\mathcal{T}}
\renewcommand{\SS}{\mathcal{S}}

\newcommand{\UC}{\mathcal{UC}}

\renewcommand{\aa}{\mathbf{a}}
\newcommand{\bb}{\mathbf{b}}

\newcommand{\mm}{\mathbf{m}}
\newcommand{\vv}{\mathbf{v}}

\newcommand{\ww}{\mathbf{w}}
\newcommand{\xx}{\mathbf{x}}
\newcommand{\yy}{\mathbf{y}}

\newcommand{\mumu}{\boldsymbol{\mu}}

\newcommand{\AAA}{\mathscr{A}}

\newcommand{\Dy}{\mathbf{Dyadic}}

\newcommand{\diam}[1]{\text{diam}\!\left({#1}\right)}
\newcommand{\supp}[1]{\mathrm{supp}\!\left({#1}\right)}
\newcommand{\absolute}[1] {\left|{#1}\right|}

\newcommand{\norm}[1]{\left\|{#1}\right\|}

\newcommand{\vol}{\mathrm{vol}}

\newcommand{\area}{\mathcal{A}}

\newcommand{\ignore}[1]  {}

\theoremstyle{plain}
\newtheorem{thm}{Theorem}[section]

\newtheorem{lem}[thm]{Lemma}

\newtheorem{ques}[thm]{Question}

\theoremstyle{definition}
\newtheorem{definition}[thm]{Definition}

\newtheorem{remark}[thm]{Remark}

\numberwithin{equation}{section}
\swapnumbers

\newif\ifdraft\drafttrue
\draftfalse
\usepackage{color}

\begin{document}
		\maketitle
\begin{abstract}
      We extend a discrepancy bound of Lagarias and Pleasants for local weight distributions on linearly repetitive Delone sets and show that a similar bound holds also for the more general case of Delone sets without finite local complexity if linear repetitivity is replaced by $\varepsilon$-linear repetitivity. 
       As a result we establish that Delone sets that are $\varepsilon$-linear repetitive for some sufficiently small $\varepsilon$ are rectifiable, and that incommensurable multiscale substitution tilings are never almost linearly repetitive.
\end{abstract}

\section{Introduction}\label{sec:introduction}
The property of linear repetitivity plays a central role in the study of mathematical models of quasicrystals and in particular in the study of aperiodic tilings and Delone sets. This is due both to the various dynamical and geometric implications of linear repetitivity, as well as to the fact that several well-studied constructions in aperiodic order are known to have this property, including  primitive self-similar tilings of finite local complexity \cite{Solomyak-recognizability} and certain cut-and-project sets \cite{Haynes-Koivusalo-Walton,Koivusalo-Walton}. In this paper we consider Delone sets and tilings of infinite local complexity, which have seen a surge of interest in recent years with examples including \cite{Danzer,Frank,Frank-Robinson,Frank-Sadun1,Frank-Sadun2,Frettloh-Richard,Lee-Solomyak,Generalized_Pinwheel,Smilansky-Solomon1} and \cite{Smilansky-Solomon2}, and for which
a suitable extension of the notion of linear repetitivity is required. Our study is motivated by the question of rectifiability of Delone sets of infinite local complexity, and in particular those defined by multiscale substitution tilings.

A set $\Lambda \subset \R^d$ is \emph{Delone} if it is \emph{uniformly discrete} and \emph{relatively dense}, that is, if there exist constants $r,R>0$ so that every ball of radius $r$ contains at most one point of $\Lambda$ and $\Lambda$ intersects every ball of radius $R$. In our setup balls and distances are taken with respect to the metric induced by the sup-norm, denoted by $\norm{\cdot}$. We also define $r_\Lambda$ and $R_\Lambda$, the \emph{packing constant} and the \emph{covering constant} of $\Lambda$, respectively, by
\[
r_\Lambda = \inf\{\norm{\xx_1 - \xx_2}  \:|\: \xx_1\neq \xx_2, \xx_1,\xx_2\in\Lambda  \}, \quad 
R_\Lambda =  \sup\{\norm{\xx - \yy}  \:|\: \xx\in\Lambda, \yy\in\R^d \}.
\]
We will assume without loss of generality that $\Lambda$ has covering radius $R_\Lambda=1$. 
For $t>0$ and $\xx\in\Lambda$, let $B(\xx,t)$ be the ball of radius $t$ centered at $\xx$, then the set $\Lambda\cap B(\xx,t)$ is the \emph{$t$-patch of $\Lambda$ at $\xx$}. A Delone set has \emph{finite local complexity (FLC)} if for every $t>0$ its collection of $t$-patches is finite modulo translations, and it is \emph{repetitive} if for every $r>0$ there exists an $R=R(r)>0$ so that every $R$-patch of $\Lambda$ contains translated copies of every $r$-patch of $\Lambda$. Repetitivity is equivalent to the minimality of the \emph{hull} of $\Lambda$, $\X(\Lambda)$, which is the orbit closure of $\Lambda$ with respect to translations, see \cite[Prop. 5.4]{BaakeGrimm} and \cite[Theorem 3.2]{Lagarias-Pleasants} for precise statements. 
Finally, a repetitive Delone set  is \emph{linearly repetitive} if $R(r)$ can be chosen to be a linear function. 

A Delone set $\Lambda\subset\R^d$ is \emph{rectifiable} if there exists a biLipschitz bijection between $\Lambda$ and $\Z^d$. While Burago and Kleiner's well-known rectifiability condition (\cite{BK2} for $d=2$, \cite{ACG} for $d\ge 3$, included below as Theorem \ref{thm:BK_condition}) can be used to establish rectifiability for a large family of constructions, examples include certain sets constructed via substitution tilings \cite{Solomon11} and the cut-and-project method \cite{HKW}, the pursue of non-rectifiable Delone sets is known to be a difficult problem. The question whether such Delone sets exist was posed by Gromov in \cite[p. 23]{Gr}, and according to \cite{BK2} it was also posed by Furstenberg in connection with Kakutani equivalence for $\R^2$-actions. It was answered on the affirmative independently by Burago and Kleiner in \cite{BK1} and by McMullen in \cite{McMullen}. Concrete examples of such Delone sets were later provided in \cite{Cortez-Navas}, based on the construction in \cite{BK1}, see also \cite{Garber} and \cite{Magazinov}. It was shown in \cite{ACG} that linear repetitive Delone sets satisfy Burago and Kleiner's condition and are therefore rectifiable. In fact, in this case the implied biLipschitz bijection can be extended to a biLipschitz homeomorphism of $\R^d$, see \cite{Navas}, though it remains unknown if this is always the case, see \cite[Problem 2.6.1]{ADG}. It follows that non-linear repetitive Delone sets emerge as natural candidates for non-rectifiable Delone sets with a minimal hull. This remains true also when moving beyond the FLC setup, as we will see below.

While sets of infinite local complexity can never be repetitive, they can nevertheless be $\varepsilon$-repetitive. Recall that the \emph{Hausdorff distance} between two compact subsets $K_1,K_2\subset\R^d$ is defined by  
\begin{equation}\label{eq:Hausdorff-metric}
D_H\left(K_1,K_2\right) = 
\inf\left\{ \varepsilon>0 \: \Big| \: \begin{matrix}K_1\subset K_2^{(+\varepsilon)}\\
K_2\subset K_1^{(+\varepsilon)}
\end{matrix}\right\},
\end{equation}
where $A^{(+\varepsilon)} = \bigcup_{\xx\in A}B(\xx,\varepsilon)$, the $\varepsilon$ neighborhood of the set $A$. We say that $K_1$ is an \emph{$\varepsilon$-copy} of $K_2$ if it is of distance at most $\varepsilon$ of some translation of $K_2$.

\begin{definition}\label{def:linear_repetitivity}
	Let $\varepsilon>0$. A Delone set $\Lambda \subset \R^d$ is \emph{$\varepsilon$-repetitive} if for every $r>0$ there exists $R=R(r,\varepsilon) >0$ such that every $R$-patch of $\Lambda$ contains an $\varepsilon$-copy of every $r$-patch of $\Lambda$. If there exists $C_{\rm{rep}}=C_{\rm{rep}}(\Lambda,\varepsilon)$ for which this holds  for $R=C_{\rm{rep}}\cdot r$, then $\Lambda $ is \emph{$\varepsilon$-linearly repetitive}. It is \emph{almost  repetitive} if it is $\varepsilon$-repetitive for every $\varepsilon>0$, and \emph{almost linearly repetitive} if it is $\varepsilon$-linearly repetitive for every $\varepsilon>0$. 
\end{definition}

Lagarias and Pleasants showed in \cite{Lagarias-Pleasants} that linear repetitivity implies certain discrepancy bounds, which in turn imply uniform patch frequency and hence unique ergodicity of the hull, see also \cite[Corollary 4.6]{Damanik-Lenz}. Further connections between linear repetitivity and dynamics  include \cite{AlistePrieto-Coronel,Besbes,Besbes-Boshernitzan-Lenz,Cortez-Durand-Petite,Damanik-Lenz} and \cite{Du}, and we refer the reader to \cite{ACCDP} for a comprehensive discussion and for many additional references. In parallel with the FLC setup,  almost repetitivity of $\Lambda$ is equivalent to the minimality of the hull of $\Lambda$ in the infinite local complexity case, see \cite[Theorem 3.11]{Frettloh-Richard}. 
Adapting ideas from \cite{Damanik-Lenz} and \cite{Lagarias-Pleasants}, and in particular the approach of \cite{Damanik-Lenz} and their weight function, Frettl\"oh and Richard showed in \cite{Frettloh-Richard} that if $\Lambda$ is almost linearly repetitive then the hull of $\Lambda$ is uniquely ergodic. 
Our first result is the following upper bound on the discrepancy for $\varepsilon$-linearly repetitive Delone sets, adapting the approach of Lagarias and Pleasants from \cite{Lagarias-Pleasants} and the study of their weight distribution functions in the FLC case.

A \emph{box} $B$ in $\R^d$ is a set of the form $\bigtimes_{i=1}^d[a_i,b_i]$ with $a_i<b_i$ for every $i$. Denote by $\vol(B)$ the Lebesgue measure of a box $B$ and by $\ell(B)=\min\{b_i-a_i \mid 1\le i\le d\}$ its width.

\begin{thm}\label{thm:Discrepancy_bound_for_eps_linear_repet.}
	Let $\Lambda \subset \R^d$ be a Delone set and assume that $\Lambda$ is $\varepsilon$-linearly repetitive for some fixed $\varepsilon<r_\Lambda$. Then there exist an asymptotic density $\mu$ and constants $\alpha>0$ and $0<\delta<1/2$ such that for every box $B$ we have
	\begin{equation}\label{eq:Discrepancy_bound_for_eps_linear_repet.}
	\absolute{\#(\Lambda \cap B) - \mu\cdot \vol(B)} \le \alpha\cdot \frac{\vol(B)}{\ell(B)^{\delta}}.
	\end{equation}
	where $\alpha$ and $\delta$ depend on $d$, $\varepsilon$ and $\Lambda$.
\end{thm}

Note that since $\delta<1/2$ the discrepancy bound in Theorem \ref{thm:Discrepancy_bound_for_eps_linear_repet.} is greater than the surface area of $B$, which for a cube would correspond to the value $\delta=1$.
An extension of this result to finite unions of unit cubes is given below in Theorem \ref{thm:discrepancy_for_UC}. These discrepancy bounds are applied to deduce the two corollaries stated below as Theorems \ref{cor:extending_BL_result_of_ACG} and \ref{thm:Multiscale_are_not_almost_linearly_repetitive}. The first is an immediate consequence of 
Theorem \ref{thm:Discrepancy_bound_for_eps_linear_repet.}, when combined with Burago and Kleiner's condition and with the argument of Navas in \cite{Navas}.

\begin{thm}\label{cor:extending_BL_result_of_ACG}
			Let $\Lambda \subset \R^d$ be a Delone set and assume that $\Lambda$ is $\varepsilon$-linearly repetitive for some fixed $\varepsilon<r_\Lambda$. Then $\Lambda$ is rectifiable and, moreover, there exists a biLipschitz homeomorphism $F:\R^d\to\R^d$ such that $F(\Lambda) = \Z^d$.
\end{thm}

Theorem \ref{cor:extending_BL_result_of_ACG} extends the results established in  \cite{ACG} and \cite{Navas} on the rectifiability of linear repetitive Delone sets. Indeed, rectifiability is a large-scale property and should not be affected by moving each point a small distance. 
	
	\begin{remark}
		A related large-scale property is uniform spreadness, where a Delone set $\Lambda\subset\R^d$ is  \emph{uniformly spread} if there exists a bijection between $\Lambda$ and $\Z^d$ that moves every point in $\Lambda$ a bounded distance. It is not hard to see that a uniformly spread Delone set is always rectifiable. We note that uniformly spreadness and $\varepsilon$-linear repetitivity imply distinct but similar discrepancy bounds, compare Laczkovich's condition \cite{Laczk} and Theorem \ref{thm:Discrepancy_bound_for_eps_linear_repet.}. The two properties are independent: the examples of non-uniformly spread Delone sets associated with primitive substitution tilings that were described in \cite{Solomon14} are linear repetitive by \cite{Solomyak-recognizability}, and a uniformly spread Delone set $\Lambda$ that is not $\varepsilon$-linear repetitive for any $\varepsilon<r_\Lambda$ can be constructed by appropriately perturbing all points in $\Z^d$.
	\end{remark}

Theorem \ref{thm:Multiscale_are_not_almost_linearly_repetitive} concerns with incommensurable multiscale substitution tilings, a class of tilings of infinite local complexity that was recently introduced by the authors in \cite{Smilansky-Solomon1}.  The construction of such tilings will be recalled in \S \ref{sec:multiscale tilings} together with the terms used in the following statement.

\begin{thm}\label{thm:Multiscale_are_not_almost_linearly_repetitive}
	Let $\sigma$ be an irreducible incommensurable multiscale substitution scheme in $\R^d$, with polytope prototiles, and let  $\TT\in\X_\sigma$ be a multiscale substitution tiling. Then $\TT$  is not $\varepsilon$-linearly repetitive for any sufficiently small $\varepsilon$. In particular, incommensurable multiscale substitution tilings are never almost linearly repetitive.
\end{thm}
Theorem \ref{thm:Multiscale_are_not_almost_linearly_repetitive} is another manner in which incommensurable multiscale substitution tilings defer from the classical construction of substitution tilings, compare with \cite{Solomyak-recognizability}. 
In view of Theorem \ref{cor:extending_BL_result_of_ACG}, Theorem \ref{thm:Multiscale_are_not_almost_linearly_repetitive} strengthens the candidacy of sets associated with incommensurable multiscale tilings for non-rectifiability. These are never uniformly spread, and furthermore, certain examples do not satisfy Burago and Kleiner's aforementioned rectifiability condition, which is already not a trivial result, see \cite[\S 8]{Smilansky-Solomon1}. While the rectifiability of primitive substitution tilings was established in \cite{Solomon11}, the following question remains:  
\begin{ques}
Does there exist a non-rectifiable multiscale substitution tiling?
\end{ques}


\begin{remark}
	As mentioned above, almost linear repetitivity of a Delone set implies the unique ergodicity of its hull, and by \cite{Smilansky-Solomon1} unique ergodicity holds also in the case of incommensurable multiscale substitution tilings. In view of this we note that Theorem \ref{thm:Multiscale_are_not_almost_linearly_repetitive} points to a large family of constructions whose hulls are uniquely ergodic despite being not almost linearly repetitive, adding to the earlier examples that appeared in  \cite{Cortez-Navas}.
\end{remark}

We note that it is our belief that the assumption on the tiles in Theorem \ref{thm:Multiscale_are_not_almost_linearly_repetitive} is not restrictive. An example of an incommensurable multiscale substitution tiling with non-polytope tiles has yet to be discovered. We believe that this is impossible, and we pose the following question:  
\begin{ques}
	Does there exist an irreducible incommensurable multiscale substitution scheme on a finite set of non-polytope prototiles?
\end{ques}

\subsection*{Acknowledgments} We thank the anonymous
referees for a careful reading of the paper and for helpful
comments and suggestions.

\section{Lagarias-Pleasants for almost linear repetitivity}\label{sec:LP-analogue}
We denote the \emph{volume} and the \emph{surface area} of a box $B=\bigtimes_{i=1}^d[a_i,b_i]$ by $\vol(B)$ and $ \area (B)$. Setting $\ell_i=b_i-a_i$ for $i=1,\ldots,d$ we have therefore
\begin{equation}\label{eq:volume and width in terms of lengths}
\vol(B)=\prod_{i=1}^d\ell_i \,,\quad  \area (B)=2\vol(B)\sum_{i=1}^d\frac{1}{\ell_i} .
\end{equation}
In addition, we define the \emph{width} and the \emph{middle point} of $B$ by  
\[  \ell(B)=\min_{1\le i\le d} \ell_i\,,\quad \mm(B) = \frac{\aa+\bb}{2}, \]
where $\aa = (a_1,\ldots,a_d)$ and $\bb = (b_1,\ldots,b_d)$. 
 Definition \ref{def:weight_distribution} below is similar to \cite[Definition 5.2]{Lagarias-Pleasants}, the only difference being the appearance of $\varepsilon>0$ in item (b).

\begin{definition}\label{def:weight_distribution}
	Let $\Lambda \subset \R^d$ be a Delone set, let $p\in\N$ and fix $\varepsilon, t_0>0$. An
	\emph{$\varepsilon$-weight distribution} is a function $\ww$ with values in $\R^p$, whose domain is the collection of all boxes $B$ in $\R^d$ with $\ell(B)\ge t_0$, for which there is a constant $C_w\ge 1$ so that the following three properties hold for every box $B$ in its domain:
	\begin{itemize}
		\item[(a)] ({\it boundedness})
		$\norm{\frac{\ww(B)}{\vol(B)}} \le C_w$.
		\item[(b)]  ({\it almost approximate invariance})
		For any $\vv\in\R^d$, if $(B\cap\Lambda)-\vv$ and $(B-\vv)\cap\Lambda$ are $\varepsilon$-copies
		then $\norm{\ww(B-\vv) - \ww(B)} \le C_w\cdot  \area (B)$. 
		\item[(c)] ({\it approximate additivity})
		If $B = B_1 \cup \cdots \cup B_k$ is a union of boxes with pairwise disjoint interiors in the domain of $\ww$, then 
		\[ \norm{\ww(B) - \sum_{j=1}^k\ww(B_j)} \le  C_w\cdot \left(\sum_{j=1}^k  \area (B_j) \right).\]
	\end{itemize}
\end{definition}

For $t>0$, we denote by $\BB(t)$ the collection of all ``squarish'' boxes in $\R^d$, boxes for which all side lengths are between $t$ and $2t$. For $\varepsilon,t_0>0$ and a Delone set $\Lambda\subset\R^d$, let $w$ be a real valued $\varepsilon$-weight distribution. Then for any $t\ge t_0$, the \emph{upper density} and \emph{lower density} of $w$ are defined by 
\begin{equation}\label{eq:upper/lower_density}
{\bf \mu}^+(t):= \sup_{B\in \BB(t)} \frac{w(B)}{\vol(B)}
\quad\text{and}\quad
{\bf \mu}^-(t):= \inf_{B\in \BB(t)} \frac{w(B)}{\vol(B)}.
\end{equation}

\begin{thm}\label{thm:LP_for_almost_linear_repetitivity}
	Let $\varepsilon>0$ and let $\Lambda \subset \R^d$ be an $\varepsilon$-linearly repetitive Delone set. Then there exist constants $alpha>0$ and $0<\delta = \delta(\Lambda,\varepsilon) <1/2$ such that every $\varepsilon$-weight distribution $\ww$  on $\Lambda$ has an asymptotic weight density $\mumu \in \R^p$ for which every box $B$ with $\ell(B)\ge 2C_{\rm{rep}}$ in the domain of $\ww$ satisfies
	\begin{equation}\label{eq:LP_for_almost_linear_repetitivity}
	\norm{\frac{\ww(B)}{\vol(B)} - \mumu} \le \alpha\cdot\ell(B)^{-\delta},
	\end{equation}
	where $\alpha$ may depend on $d$, $\varepsilon$, $\Lambda$ and $\ww$, and $C_{\rm{rep}} \ge1$ is as in Definition \ref{def:linear_repetitivity}.
\end{thm}

Theorem \ref{thm:LP_for_almost_linear_repetitivity} is analogous to \cite[Theorem 5.1]{Lagarias-Pleasants}. The proof can be extended to our more general context, but since some small changes are needed we include it below. The proof requires the following Lemma \ref{lem:LP_Lemma5.1}, which was established by Lagarias and Pleasants relying only on properties (a) and (c) of local weight distributions in \cite[Definition 5.2]{Lagarias-Pleasants}. Since these are identical to properties (a) and (c) in our Definition  \ref{def:weight_distribution} of $\varepsilon$-weight distributions given above, we do not repeat the proof.
\begin{lem}[\cite{Lagarias-Pleasants}, Lemma 5.1]\label{lem:LP_Lemma5.1}
	Let $\varepsilon>0$ and let $w$ be a real-valued $\varepsilon$-weight distribution on a Delone set $\Lambda \subset \R^d$. Then the limits $\mu^+ = \lim_{t\to\infty} \mu^+(t)$ and $\mu^- = \lim_{t\to\infty} \mu^-(t)$ exist. 
\end{lem}

\begin{proof}[Proof of Theorem \ref{thm:LP_for_almost_linear_repetitivity}]
	Let $\ww$ be an $\varepsilon$-weight distribution defined on all boxes with $\ell(B)\ge t_0$.  Considering each coordinate individually, it suffices to prove the assertion under the assumption that $\ww = w$ is real-valued. 
	
	First, we show that $\mu^+=\mu^-$. Let $t\ge t_1:=\max\{t_0,2C_{\rm{rep}}\}$, then by the definition of $\mu^-(t)$ there exists some box $B_1'\in \BB(t)$ that satisfies 
	\begin{equation}\label{eq:LP-analogue1}
	\frac{w(B_1')}{\vol(B_1')} \le \mu^-(t) + \frac{1}{t}.
	\end{equation}
	Set $B_2' := B(\mm(B_1'),2t)$, a ball of radius $2t$ (in the sup-norm) that contains $B_1'$ positioned such that $B_1'$ is at distance of at least $t$ from the boundary of $B_2'$. By $\varepsilon$-linear repetitivity, every ball of radius $C_{\rm{rep}} \cdot 2t$, and in particular every box $B\in\BB(C_{\rm{rep}} \cdot 2t)$, contains a ball $B_2$ of radius $2t$ for which $B_2\cap\Lambda$ is an $\varepsilon$-copy of $B'_2\cap\Lambda$.
	In particular, if we set $B_1:=B_1' - \mm(B_2') + \mm(B_2) \subset B_2$, then $B_1\cap\Lambda$ is an $\varepsilon$-copy of $B'_1\cap\Lambda$.
	By property (b) of Definition \ref{def:weight_distribution} with $\vv = - \mm(B_2') + \mm(B_2)$ we obtain 
	\begin{equation}\label{eq:LP-analogue2}
	\absolute{w(B_1) - w(B_1')} \le C_w\cdot  \area (B_1). 
	\end{equation}
	
	Combining \eqref{eq:LP-analogue1} and \eqref{eq:LP-analogue2}, every box $B\in\BB(C_{\rm{rep}} \cdot 2t)$ contains a box $B_1\in \BB(t)$ that is positioned inside $B$ at distance of at least $t$ from its boundary, that satisfies
	\begin{equation}\label{eq:LP-analogue3}
	\frac{w(B_1)}{\vol(B_1)} \le \mu^-(t) + \frac{2dC_w+1}{t},
	\end{equation}
	where the upper bound $2d/t$ of $ \area (B_1)/\vol(B_1)$ follows from \eqref{eq:volume and width in terms of lengths}. Consider a partition $B=B_1\cup\ldots\cup B_k$ of $B$ into boxes with pairwise disjoint interiors in $\BB(t)$, one of which is $B_1$. Such a partition exists because of the way $B_1$ is positioned inside $B$. Denote $C_1 = \vol(B_1)/\vol(B)$. Combining \eqref{eq:LP-analogue3} and property (c) of Definition \ref{def:weight_distribution} we get
\begin{equation}\label{eq:LP-analogue4}
\begin{aligned}
\frac{w(B)}{\vol(B)} &\le \frac{\sum_{j=1}^kw(B_j) + C_w\left(\sum_{j=1}^k  \area (B_j)\right)}{\vol(B)} \\
&\le 
C_1\frac{w(B_1)}{\vol(B_1)} + \sum_{j=2}^k\left(\frac{w(B_j)}{\vol(B_j)}\frac{\vol(B_j)}{\vol(B)}\right) +C_w\left(\sum_{j=1}^k \frac{ \area (B_j)}{\vol(B_j)}\frac{\vol(B_j)}{\vol(B)}\right) \\
&\le C_1\left(\mu^-(t)+\frac{2dC_w+1}{t}\right) +(1-C_1)\mu^+(t)+\frac{2dC_w}{t}.
\end{aligned}
\end{equation}
This inequality holds for all $B\in\BB(C_{\rm{rep}} \cdot 2t)$, thus 
\begin{equation}\label{eq:LP-analogue5}
\mu^+(C_{\rm{rep}} \cdot 2t)\le C_1\left(\mu^-(t)+\frac{2dC_w+1}{t}\right) + (1-C_1)\mu^+(t)+\frac{2dC_w}{t}.
\end{equation} 
Letting $t\to\infty$ and in view of Lemma \ref{lem:LP_Lemma5.1}, we conclude that $\mu^+ \le \mu^-$ and so $\mu:=\mu^+=\mu^-$.  
	
	We now bound the error term. Repeating the above argument, switching the rolls of $\mu^+$ and $\mu^-$ and replacing $B_1$ and $C_1$ with suitable $\widetilde{B}_1$ and $\widetilde{C}_1$ yields 
	\begin{equation}\label{eq:LP-analogue6}
	\mu^-(C_{\rm{rep}} \cdot 2t)\ge \widetilde{C}_1\left(\mu^+(t)-\frac{2dC_w+1}{t}\right) + (1-\widetilde{C}_1)\mu^-(t)-\frac{2dC_w}{t}.
	\end{equation}
	Denote $C_2=2C_{\rm{rep}}$ and define $\Delta(t):=\mu^+(t) - \mu^-(t)$. Since $(2C_2)^{-d} \le C_1,\widetilde{C}_1\le (C_2/2)^{-d}$,  \eqref{eq:LP-analogue5} and \eqref{eq:LP-analogue6} imply that 
	\begin{equation}\label{eq: LP 5.6}
	\Delta(C_2\cdot t)\le (1-C_1-\widetilde{C}_1)\Delta(t) + \frac{8dC_w}{t} \le (1-2\cdot (2C_2)^{-d})\Delta(t) + \frac{8dC_w}{t},
	\end{equation}
	for every $t\ge t_1$. 
	
	Next, we apply the relation on $\Delta(t)$ in \eqref{eq: LP 5.6} to establish the result for the case $B\in\BB(t)$ with $t\ge t_1$. Let $C_3>0$ be a constant chosen so that
	\begin{equation}\label{eq:LP-analogue7}
	\Delta(t)\le C_3\cdot t^{-\delta}
	\end{equation}
	holds for  every $t_1\le t\le C_2\cdot t_1$ and so that  $C_3 >8dC_w \cdot (2C_2)^{d} $, where 
	\begin{equation}\label{delta}
		\delta = \frac{\log(1/(1-(2C_2)^{-d}))}{\log C_2}.
	\end{equation}
	It is straightforward to check that since $C_{\rm{rep}},d\ge1$ we have $\delta\le\log(4/3)/\log2<1/2$. Assume as an induction hypothesis that \eqref{eq:LP-analogue7} holds for all $t_1\le t\le C_2^kt_1$, which holds for $k=1$, and let $t_1\le t\le C_2^{k+1}t_1$. Then by \eqref{eq: LP 5.6} 
	\begin{equation*}
	\Delta(t)\le(1-2(2C_2)^{-d})\Delta(t/C_2)+\frac{8dC_w}{t/C_2}.
	\end{equation*}
	Note that since $C_{\rm{rep}}\ge1$ and $C_2=2C_{\rm{rep}}$ we deduce that $0<\delta<1$. Combined with the induction hypothesis on $C_3$, and since $t\ge C_2$, the inequality can be extended to get
		\begin{equation*}
	\Delta(t)\le(1-(2C_2)^{-d})C_3(t/C_2)^{-\delta}=C_3t^{-\delta},
	\end{equation*}
	and so by induction \eqref{eq:LP-analogue7}  holds for all $t\ge t_1$. 
	
	Note that for every $s\ge t\ge t_1$, every box $B\in\BB(s)$ can be subdivided into boxes of $\BB(t)$. By a  computation similar to \eqref{eq:LP-analogue4} one sees that $\mu^+(s)\le \mu^+(t)+d2^{d}C_w/t$ and $\mu^-(s)\ge \mu^-(t)-d2^{d}C_w/t$, where we use property (c) of Definition \ref{def:weight_distribution} and naive bounds on the volumes and surface areas of boxes in $\BB(s)$ and $\BB(t)$. 
	Since this is true for arbitrarily large values of $s$ we deduce that for $C_4=d2^{d}C_w$
	\begin{equation*}
		\mu^-(t)-\frac{C_4}{t} \le \mu \le \mu^+(t)+\frac{C_4}{t}. 
	\end{equation*}
	In view of \eqref{eq:LP-analogue7}, for every box $B\in\BB(t)$ for $t\ge t_1$ we have obtained
	\[\absolute{\frac{w(B)}{\vol(B)} - \mu}\le \Delta(t) + \frac{C_4}{t}\le C_5\cdot t^{-\delta}\]
	for a constant $C_5$ that depends on $w$, $\varepsilon$, $\Lambda$ and $d$, implying the assertion for $B\in\BB(t)$.
	
	For an arbitrary box $B$ with width $\ell(B)\ge t_1$, we partition $B$ into boxes $B_i$ in $\BB(\ell(B))$. Then by property (c) of Definition \ref{def:weight_distribution}, and using simple approximations similar to those mentioned above, there exists $\alpha>0$, that depends on $w$, $\varepsilon$, $\Lambda$ and $d$, so that for any $B$ with width $\ell(B)\ge t_1$ we have
	\begin{equation*}
	\absolute{\frac{w(B)}{\vol(B)} - \mu}\le \absolute{\frac{w(B)}{\vol(B)} - \frac{\sum_iw(B_i)}{\vol(B)}} + \absolute{\frac{\sum_iw(B_i)} {\vol(B)}- \mu}\le \frac{C_wd2^{d}}{\ell(B)}+\frac{C_5}{\ell(B)^{\delta}}\le\alpha\cdot\ell(B)^{-\delta},
	\end{equation*}
	finishing the proof.
\end{proof}

In view of the above, our main result follows. Given a Delone set $\Lambda\subset\R^d$ we consider the function 
\begin{equation}\label{eq:kappa}
N_{\Lambda}(A):=\#(A\cap\Lambda).
\end{equation}

\begin{proof}[Proof of Theorem \ref{thm:Discrepancy_bound_for_eps_linear_repet.}]
	Observe that $N_{\Lambda}(A)$ is defined on every subset $A\subset\R^d$, and its restriction to boxes is clearly an $\varepsilon$-weight distribution for any $0<\varepsilon<r_\Lambda$. It is enough to prove the assertion for boxes $B$ with sufficiently large width $\ell(B)$, which follows directly from Theorem \ref{thm:LP_for_almost_linear_repetitivity} upon multiplying by $\vol(B)$.
\end{proof}

\section{Rectifiability of $\varepsilon$-linearly repetitive, non-FLC, Delone sets}
The following is an equivalent reformulation of Burago and Kleiner's sufficient condition for rectifiability, established for $d=2$ in \cite{BK2}  and for $d\ge 2$ in \cite{ACG}. 

\begin{thm}[\cite{ACG}, \cite{BK2}]\label{thm:BK_condition}
Let $\Lambda\subset \R^d$ be a Delone set. If there exists $\rho>0$ for which the sum
\begin{equation}\label{eq:BK_condition}
\sum_{k=1}^\infty \left[ \sup_{\xx\in\Z^d}\frac{\absolute{\#(\Lambda\cap B(\xx,2^k))-\rho \cdot \vol\left( B(\xx,2^k)\right)}}{\vol\left( B(\xx,2^k)\right)}\right] 
\end{equation}
is convergent, then $\Lambda$ is rectifiable.
\end{thm}

\begin{proof}[Proof of Theorem \ref{cor:extending_BL_result_of_ACG}]
	Let $\Lambda \subset \R^d$ be an $\varepsilon$-linearly repetitive Delone set with $\varepsilon<r_\Lambda$. Applying Theorem \ref{thm:Discrepancy_bound_for_eps_linear_repet.} and plugging $\rho=\mu$ in \eqref{eq:BK_condition} yields the series 
	\[
	\alpha\sum_{k=1}^\infty 2^{-k\delta}, 
	\]
	where $\delta>0$ is fixed and depends on $\varepsilon$ and $ \Lambda$. Then by Theorem \ref{thm:BK_condition} rectifiability follows. 
	
	Observe that the second part of the statement of Theorem \ref{cor:extending_BL_result_of_ACG} is an immediate consequence as well. It was pointed out by Navas in \cite{Navas} that for every linearly repetitive Delone set $\Lambda\subset\R^d$ there is a biLipschitz homeomorphism $F:\R^d\to\R^d$ satisfying $F(\Lambda)=\Z^d$. In fact, the FLC assumption played no part in the proof, and it was actually shown that the above holds for any $\Lambda$ that satisfies Burago and Kleiner's condition in Theorem \ref{thm:BK_condition}. Then in view of the proof of the first part of Theorem \ref{cor:extending_BL_result_of_ACG}, the second part is obtained.
\end{proof}

\section{Discrepancy bounds for unions of cubes}
In this section we extend the discrepancy bound established in \S\ref{sec:LP-analogue}  to sets that are finite unions of unit lattice cubes.
We denote by 
\[
\CC_d:=\left\{ [a_1,a_1+1)\times\ldots\times[a_d,a_d+1) \:|\: (a_1,\ldots,a_d)\in\Z^d \right\} 
\]
the set of all half-closed unit lattice cubes in $\R^d$, by $\UC_d$ the collection of all finite unions of elements of $\CC_d$, and by $\vol(U)$ and $\area(U)$ the volume and surface area, respectively, of an element $U\in\UC_d$. Let
\[
\Dy_d:=\left\{ [2^ka_1,2^ka_1+2^k)\times\ldots\times[2^ka_d,2^ka_d+2^k) \:|\: k\in\Z_{\ge0}, (a_1,\ldots,a_d)\in\Z^d \right\}  
\]
denote the set of all half-closed dyadic cubes in $\R^d$ with vertices in $2^k\Z^d$, for some $k\in\N$. The following notion was introduced in \cite[p. 41]{Laczk}. 
Given a collection $\AAA$ of elements of $\Dy_d$,  define $\SS(\AAA)$ to be the closure of $\AAA$ under the
operations of disjoint union and proper difference with the restriction that each
element of $\AAA$ can be used at most once. 
We rely on the following result by Laczkovich.
\begin{lem}[\cite{Laczk}, Lemma 3.2]\label{lem:Laczkovich_Lemma3.2}
	Let 
	\begin{equation}\label{eq:Laczkovich_3.2_assumptions}
	U\in\UC_d, \quad B\in\Dy_d, \quad \text{such that } \quad U\subset B, \quad 
	\vol(U)\le \frac12 \vol(B). 
	\end{equation}
	Then there exist $B_1,\ldots, B_m \in \Dy_d$ contained in $B$, so that $U\in \SS\left(\{B_1,\ldots,B_m\}\right)$ and  
	\begin{equation}\label{eq:Laczkovich_3.2}
	\#\{i \:|\: \ell(B_i)=2^k\} \le C_6\cdot \frac{ \area (U)}{2^{k(d-1)}}
	\end{equation} 
	for every $k\in\Z_{\ge0}$, where $C_6$ depends only on the dimension $d$.
\end{lem}

Let $N_\Lambda$ be as in \eqref{eq:kappa}. The main result of this chapter is the following.
\begin{thm}\label{thm:discrepancy_for_UC}
	Let $\varepsilon>0$ and let $\Lambda \subset \R^d$ be an $\varepsilon$-linearly repetitive Delone set. Let $U\in\UC_d$ and let $B\in\Dy_d$ that relates to $U$ as in \eqref{eq:Laczkovich_3.2_assumptions}. Then
	\begin{equation}\label{eq:goal_LP_for_extended_distributions}
	\absolute{N_{\Lambda}(U) - \mu\cdot \vol(U)} \le \beta\cdot \ell(B)^{1-\delta}\cdot  \area (U),
	\end{equation}
	where $\delta$ and $\mu$ are as in Theorem \ref{thm:LP_for_almost_linear_repetitivity} and $\beta$ depends on $d$, $\varepsilon$ and $\Lambda$. 
\end{thm}

\begin{remark}
	The proof of Theorem \ref{thm:discrepancy_for_UC} holds for other $\varepsilon$-weight distributions $\ww$ that are defined on elements of $\UC_d$ in a similar way to Definition \ref{def:weight_distribution} that also satisfy
	\[\ww(U_1\cup U_2) = \ww(U_1) + \ww(U_2)\]
	for all disjoint $U_1,U_2\in\UC_d$. Another example for such a function is the patch counting function $N_{\Lambda,P}(U)$ that counts the number of centers of a given patch $P$ in the set $U$. Also note that in the particular case that $U$ is a box in $\BB(t)$ for some $t$, the bound in \eqref{eq:goal_LP_for_extended_distributions} differs from the bound given in Theorem \ref{thm:LP_for_almost_linear_repetitivity} by a constant only. 
\end{remark}

\begin{lem}\label{lem:discrepancy_on_UC1}
	For every $U\in \UC_d$, if $B_1,\ldots,B_m\in\Dy_d$ and $U\in \SS(\{B_1,\ldots,B_m  \})$ then for any $\rho\in\R$ we have 
	\begin{equation}\label{eq:discrepancy_on_UC1}
	\absolute{N_{\Lambda}(U) - \rho\cdot\vol(U)} \le \sum_{i=1}^m\absolute{N_{\Lambda}(B_i) - \rho\cdot\vol(B_i)}.
	\end{equation}
\end{lem}
\begin{proof}
	The proof is straightforward from the definition of $\SS(\{B_1,\ldots,B_m\})$ and the following two simple observations, that hold for every $\rho>0\in\R$ and every $U_1,U_2\in \UC_d$:\\
	If $U_1\cap U_2 = \varnothing$ then 
	\[
	\begin{aligned}
	\absolute{N_{\Lambda}(U_1\cup  U_2) - \rho\cdot \vol(U_1\cup  U_2)}\le \absolute{N_{\Lambda}(U_1) - \rho\cdot \vol(U_1)} + \absolute{N_{\Lambda}(U_2) - \rho\cdot \vol(U_2)},
	\end{aligned}
	\]
	and if $U_2\subset U_1$ then 
	\[
	\begin{aligned}
	\absolute{N_{\Lambda}(U_1\smallsetminus U_2) - \rho\cdot \vol(U_1\smallsetminus U_2)} \le \absolute{N_{\Lambda}(U_1) - \rho\cdot \vol(U_1)} + \absolute{N_{\Lambda}(U_2) - \rho\cdot \vol(U_2)},
	\end{aligned}
	\]
	and the result follows.
\end{proof}

\begin{proof}[Proof of Theorem \ref{thm:discrepancy_for_UC}]
	Let $U\in\UC_d$. Applying Lemma \ref{lem:Laczkovich_Lemma3.2} with respect to some cube $B$ that satisfies \eqref{eq:Laczkovich_3.2_assumptions}, we obtain $B_1,\ldots, B_m \in \Dy_d$ that are contained in $B$, and for which $U\in \SS\left(\{B_1,\ldots,B_m\}\right)$ and \eqref{eq:Laczkovich_3.2} is satisfied for every $k\in\N$. In addition, \eqref{eq:discrepancy_on_UC1} is also satisfied by Lemma \ref{lem:discrepancy_on_UC1}. Applying Theorem \ref{thm:Discrepancy_bound_for_eps_linear_repet.} on each of the boxes $B_i$ in \eqref{eq:discrepancy_on_UC1} with $\rho = \mu$, and using \eqref{eq:Laczkovich_3.2} and the formula for a geometric sum, we obtain 
	\[
	\begin{aligned}
	\absolute{N_{\Lambda}(U)-\mu\cdot\vol(U)} &\le \alpha\sum_{i=1}^m\frac{\vol(B_i)}{\ell(B_i)^\delta} \\
	&
	\le \alpha C_6 \sum_{k=0}^{\log_2\ell(B)}\frac{2^{kd}}{2^{k\delta}\cdot 2^{k(d-1)}} \area (U) \\
	&=\alpha C_6 \frac{2^{1-\delta}\ell(B)^{1-\delta}-1}{2^{1-\delta}-1} \area (U)  \\ 
	&\le\beta \cdot \ell(B)^{(1-\delta)} \area (U),
	\end{aligned}
	\]
	where $\beta>0$ depends on $d$, $\varepsilon$ and $\Lambda$. 
\end{proof}

\section{Incommensurable multiscale substitution tilings are not almost linearly repetitive}\label{sec:multiscale tilings}
A {\it multiscale substitution scheme} $\sigma$ in $\R^d$ consists of a finite set $\tau_\sigma=(T_1,\ldots,T_n)$ of  {\it prototiles} of unit volume, and {\it substitution rules} $\varrho(T_i)$ each a partition of $T_i\in\tau_\sigma$ into finitely many rescaled copies of elements of $\tau_\sigma$. 
Tilings of $\R^d$ arise as limits of patches that are generated by $\sigma$ in the following way, where a \emph{patch} is a finite union of tiles.  Position a prototile $T_i\in\tau_\sigma$ around the origin, and define the patch $F_t(T_i)$ by inflating $T_i$ by a factor of $e^t$, while substituting every tile that appears in the process according to ${\sigma}$ once its volume is greater than the unit volume. A new tile that arises as a rescaled copy of a prototile $T_j\in\tau_\sigma$ is said to be of \emph{type} $j$.
The patches $\{F_t(T_i):\,t\ge0\}$ exhaust the space, and limits taken with respect to the natural topology on the space $\mathscr{C}(\R^d)$ of closed subsets of $\R^d$, which is closely related to the Hausdorff distance described in \eqref{eq:Hausdorff-metric} and is discussed in more detail for example in \cite{Frettloh-Richard,Smilansky-Solomon1,Smilansky-Solomon2}, define a tiling space $\X_\sigma$ of multiscale substitution tilings of $\R^d$. 
A multiscale substitution scheme $\sigma$ is \emph{irreducible} if for every $1\le i,j\le n$ there exists $t>0$ so that $F_t(T_i)$ contains a tile of type $j$. It is \emph{incommensurable} if there exists a prototile $T_i\in\tau_\sigma$ and $t_1,t_2>0$ so that $t_1\notin t_2\Q$, and $F_{t_1}(T_i)$ and $F_{t_2}(T_i)$ both contain a copy of the prototile $T_i$. A patch of an irreducible incommensurable multiscale substitution tiling in $\R^2$ is illustrated below in Figure \ref{fig: patch}.		
For more details, examples, illustrations and equivalent definitions of incommensurability, multiscale substitutions schemes and the geometric objects they generate, the reader is referred to \cite{Yotam Kakutani} and \cite{Smilansky-Solomon1}.

\begin{figure}[ht!]
	\includegraphics[scale=0.9]{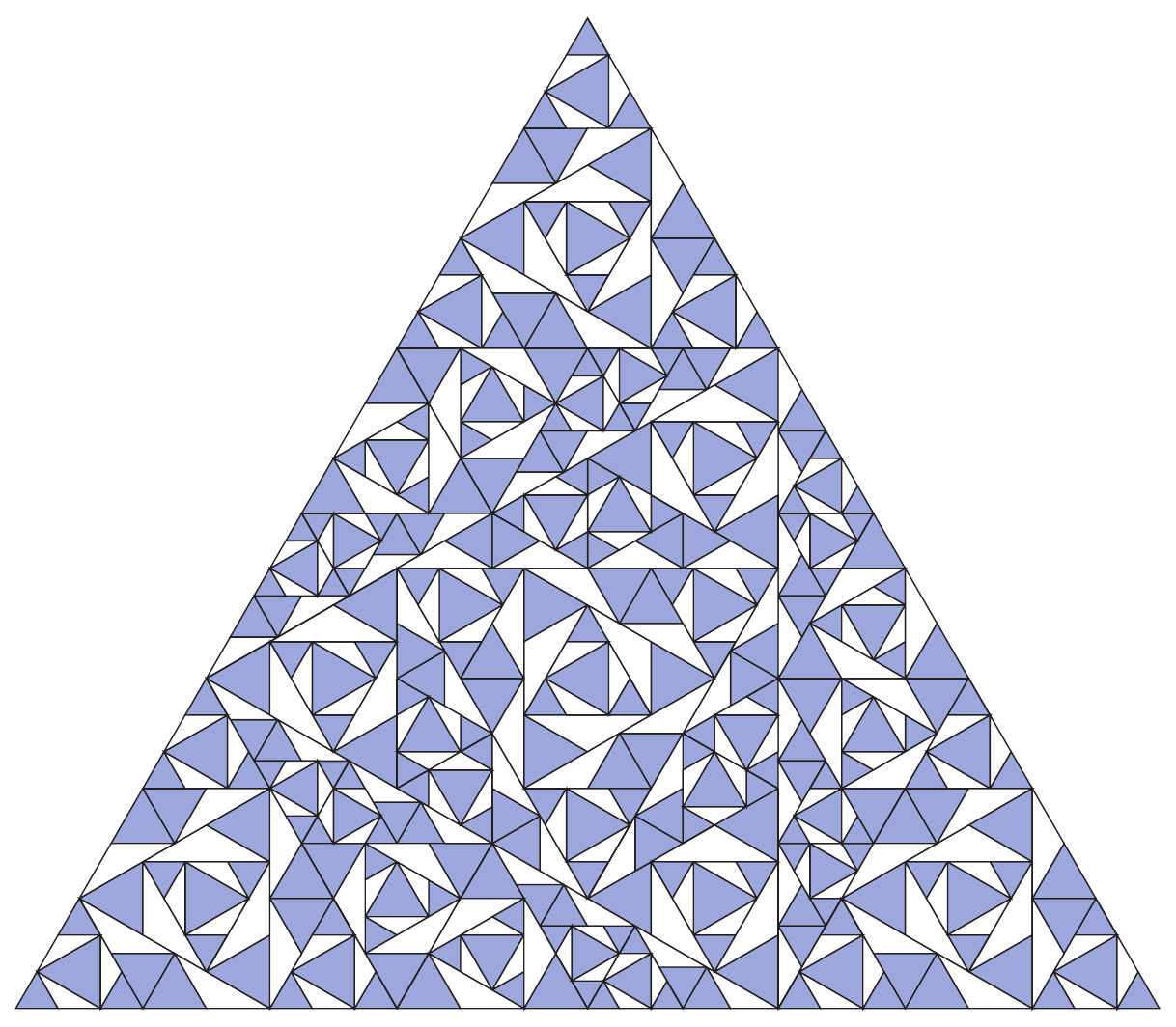}\caption{A patch of an incommensurable multiscale substitution tiling.}	\label{fig: patch}
\end{figure}

Let $\TT$ be a tiling of $\R^d$. For  $A\subset \R^d$ we denote the patch that consists of all the tiles of $\TT$ that intersect $A$ by $[A]^\TT$. Given $t>0$ and $\xx\in\R^d$, the \emph{t-patch} of $\TT$ at $\xx$ is the patch $[B(\xx,t)]^\TT$. For a patch $P$ in $\TT$ we denote by $\supp P$ the \emph{support} of $P$, which is the subset of $\R^d$ that is covered by the tiles in $P$, by $\partial P$ the union of all the boundaries of tiles in $P$ and by $\#P$ the number of tiles it consists of. We say that a patch $P_1$ is an $\varepsilon$-copy of a patch $P_2$ if $\partial P_1$ is of distance at most $\varepsilon$ of some translate of $\partial P_2$, with respect to the Hausdorff distance \eqref{eq:Hausdorff-metric}.

Suppose that every tile $T$ in a tiling $\TT$ is assigned with a type $1\le i\le n$ so that $T$ is similar to a prototile $T_i$, as is the case with multiscale substitution tilings. Marking a single point in the interior of each of the $n$ prototiles gives rise to a Delone set $\Lambda_\TT$, where each point of $\Lambda_\TT$ is contained in a distinct tile of $\TT$, with position relative to the position of the marked point in the associated prototile. More precisely, if $T=g_T(T_i)$ for a similarity $g_T$ of $\R^d$, and $\xx_i\in T_i$ is the marked point in the prototile $T_i$, then the corresponding point in $\Lambda_\TT$ is $g_T(\xx_i)\in T$. The $\varepsilon$-linear repetitivity of all Delone sets $\Lambda_\TT$ defined by the above procedure, as defined by  Definition \ref{def:linear_repetitivity}, is equivalent to the $\varepsilon$-linear repetitivity of the tiling $\TT$ as given by Definition \ref{def:tiling_eps-linear_repetitive} below.

\begin{definition}\label{def:tiling_eps-linear_repetitive}
	Let $\varepsilon>0$. A tiling $\TT$ of $\R^d$ is \emph{$\varepsilon$-linearly repetitive} if for every $r>0$ there exists $C_{\rm{rep}}=C_{\rm{rep}}(\Lambda,\varepsilon)$ such that every $(C_{\rm{rep}}\cdot r)$-patch of $\TT$ contains an $\varepsilon$-copy of every $r$-patch of $\TT$. It is \emph{almost linearly repetitive} if it is $\varepsilon$-linearly repetitive for every $\varepsilon>0$. 
\end{definition}

Consider an irreducible incommensurable multiscale substitution scheme $\sigma$ in $\R^d$. For our proof of Theorem \ref{thm:Multiscale_are_not_almost_linearly_repetitive} we will need the following two results from \cite{Smilansky-Solomon1}.

\begin{lem}[\cite{Smilansky-Solomon1}, Theorem 6.1] \label{lem:Multiscale_is_minimal}
	The dynamical system $(\X_\sigma,\R^d)$ is minimal.
\end{lem}

\begin{lem} [\cite{Smilansky-Solomon1}, Lemma 8.5]\label{lem:Multiscale_discrepancy_lower_bound}
	  For every $t_0>0$ and prototile $T_i\in\tau_\sigma$ there exist $t\ge t_0$ and $\varepsilon_0>0$ such that for every $\varepsilon\in(0,\varepsilon_0]$ we have 
	\begin{equation*}\label{eq:discrepancy_lower_bound}
	\# F_{t+\varepsilon}(T_i) - \# F_{t}(T_i) \ge C_7\cdot\frac{e^{td}}{t^k},	
	\end{equation*}	
	where $C_7>0$ and $0<k\in\N$  depend only on the parameters of $\sigma$.
\end{lem}

\begin{proof}[Proof of Theorem \ref{thm:Multiscale_are_not_almost_linearly_repetitive}] 
	In view of Lemma \ref{lem:Multiscale_is_minimal} it is enough to show that there exists a tiling $\TT\in\X_\sigma$ that is not $\varepsilon$-linearly repetitive for any sufficiently small $\varepsilon$.

	Pick a prototile $T\in\tau_\sigma$, and consider a patch of the form $F_t(T)$ for some $t>0$. Note that for every patch of the form $F_t(T)$ there exists some $\TT\in\X_\sigma$ that contains a translated copy of $F_t(T)$, positioned so that $\supp{F_t(T)}$ covers the origin (see \cite[equation (4.5)]{Smilansky-Solomon1}). Let $\TT\in\X_\sigma$ be such a tiling and let $\Lambda_\TT$ be a Delone set associated with $\TT$ in the way described above. We have 
	\[\vol(\supp{ F_t(T)} ) = \vol(T)\cdot e^{td}=e^{td}, \quad \text{and}\quad \diam{F_t(T)} = \diam{T}\cdot e^{t}, \]
	where $\diam{P}$ is the diameter of the support of the patch $P$. Note that since $T$ is a polytope, the boundary of  $\supp{F_t(T)}$ has finite $(d-1)$-dimensional Lebesgue measure, and therefore also finite $(d-1)$-dimensional Hausdorff measure. Let $U_t\in\UC_d$ denote the union of all unit lattice cubes that intersect $\supp{F_t(T)}$. By standard Hausdorff measure arguments, see e.g. \cite[p. 57]{Mattila}, we have
	\begin{equation}\label{eq:multiscale_not_linear_repetitive1}
	\area (U_t)\le C_8\cdot e^{t(d-1)},
	\end{equation}
	where $C_8$ depends on $d$ and $\sigma$. Let $B'_t\in \Dy_d$ be the smallest dyadic cube that contains $U_t$ and let $B_t\in \Dy_d$ be the dyadic cube with $\ell(B_t) = 2\ell(B'_t)$ that contains $B'_t$, then $U_t$ and $B_t$ satisfy  the requirements  \eqref{eq:Laczkovich_3.2_assumptions}. In addition, note that  
	\begin{equation}\label{eq:multiscale_not_linear_repetitive2}
	2\diam{T}\cdot e^t \le \ell(B_t)\le C_9\cdot e^t,
	\end{equation}
	where $C_9$ depends on the parameters of $\sigma$.  
	
	Let $\varepsilon<r_{\Lambda_\TT}$ and assume by way of contradiction that $\Lambda_\TT$ is $\varepsilon$-linearly repetitive. Applying Theorem  \ref{thm:discrepancy_for_UC} with $U_t$ and $B_t$, combined with \eqref{eq:multiscale_not_linear_repetitive1} and \eqref{eq:multiscale_not_linear_repetitive2}, we obtain
	\begin{equation}\label{eq:multiscale_not_linear_repetitive3}
	\absolute{N_{\Lambda_\TT}(U_t) - \mu\cdot \vol(U_t)} \le \beta\cdot \ell(B_t)^{1-\delta}\cdot \area (U_t)\le\ C_{10}\cdot e^{t(d-\delta)},
	\end{equation}
	where $C_{10}$ depends on $d$, $\varepsilon$ and $\sigma$, and $N_{\Lambda_\TT}(A) = \#(A\cap\Lambda_\TT)$ for $A\subset\R^d$.  Both $\absolute{N_{\Lambda_\TT}(U_t) - \#F_t(T)}$ and $\absolute{\vol(U_t) - \vol(\supp{F_t(T)})}$ are bounded by a constant times $ \area (U_t)$, and so we deduce that there is a constant $C_{11}$ for which 
	\begin{equation}\label{eq: upper bound}
	\absolute{ \#F_t(T) - \mu\cdot\vol(\supp{F_t(T)})} \le C_{11}\cdot e^{t(d-\delta)},	
	\end{equation}
	and this holds for any arbitrarily large $t>0$.
	
	On the other hand, combining Lemma \ref{lem:Multiscale_discrepancy_lower_bound} and the triangle inequality, and since $T$ is a polytope, for every $t_0>0$ there exists $t\ge t_0$ for which 
	\begin{equation*}\label{eq:discrepancy_lower_bound2}
	\absolute{\# F_{t}(T) - \mu\cdot \vol(\supp{ F_{t}(T)} )} \ge C_{12}\cdot\frac{e^{td}}{t^k},	
	\end{equation*}	
	where $C_{12}$ and $k\in\N_{\ge 1}$  depend only on $\sigma$. This contradicts \eqref{eq: upper bound}, completing the proof.

\end{proof}

\begin{remark}
For explicit formulas for the implied asymptotic density $\mu$ in terms of $\sigma$ see \cite{Frequencies}. We note that one may also consider complexity estimates similar to those that appear in \cite[\S5,\S A.5]{Frank-Sadun2} to show that the number of distinct patches in an incommensurable multiscale substitution tiling $\TT$ that can appear inside a big ball, up to distance $\varepsilon$, is of order strictly greater than the volume of the ball. This is of course impossible if  $\TT$ is $\varepsilon$-linear repetitive, thus offering another approach to Theorem \ref{thm:Multiscale_are_not_almost_linearly_repetitive}.
\end{remark}


\begin{thebibliography}{99}
\bibitem[ADG$^+$]{ADG} F. Adiceam, D. Damanik, F. G\"ahler, U. Grimm, A. Haynes, A. Julien, A. Navas, L. Sadun, B. Weiss, \emph{Open problems and conjectures related to the theory of mathematical quasicrystals}, Arnold Math. J. 2(4), 579--592, (2016).


\bibitem[AC]{AlistePrieto-Coronel} J. Aliste-Prieto, D. Coronel, \emph{Tower systems for linearly repetitive Delone sets}, Ergo. Theo. Dynam. Sys. 31(6), 1595--1618, (2011).

\bibitem[ACCDP]{ACCDP} J.  Aliste-Prieto, D. Coronel, M. I. Cortez,
F. Durand, S. Petite, \emph{Linearly repetitive Delone sets}, in \underline{Mathematics of aperiodic order}, eds. J. Kellendonk, D. Lenz, J. Savinien, Progr. Math. 309, Birkh\"auser/Springer, Basel, 195--222, (2015).

\bibitem[ACG]{ACG} J. Aliste-Prieto, D. Coronel, J. M. Gambaudo, \emph{Linearly repetitive Delone sets are rectifiable}, Ann. Inst. H. Poincar\'{e} Anal. Non Lin\'{e}aire 30(2), 275--290, (2013). 


\bibitem[BG]{BaakeGrimm} M. Baake, U. Grimm, \underline{Aperiodic order. Volume 1: A mathematical invitation}, Cambridge University Press, Cambridge, (2013).

\bibitem[B]{Besbes} A. Besbes, \emph{Uniform ergodic theorems on aperiodic linearly repetitive tilings and applications}, Rev. Math. Phys. 20(5), 597--623, (2008).

\bibitem[BBL]{Besbes-Boshernitzan-Lenz} A. Besbes, M. Boshernitzan, D. Lenz, \emph{Delone sets with finite local complexity: linear repetitivity versus positivity of weights}, Disc. Comp. Geom. 49(2), 335--347, (2013).
 
\bibitem[BK1]{BK1} D. Burago, B. Kleiner, \emph{Separated nets in Euclidean space and Jacobians of biLipschitz maps}, Geom. Func. Anal. 8,(2), 273--282, (1998).

\bibitem[BK2]{BK2} D. Burago, B. Kleiner, \emph{Rectifying separated nets}, Geom. Func. Anal. 12, 80--92, (2002).


\bibitem[CDP]{Cortez-Durand-Petite} M. I. Cortez, F. Durand, S. Petite, \emph{Linearly repetitive Delone systems have a finite number of nonperiodic Delone system factors}, Proc. Amer. Math. Soc. 138(3), 1033--1046, (2010).

\bibitem[CN]{Cortez-Navas} M. I. Cortez, A. Navas, \emph{Some examples of repetitive, non-rectifiable Delone sets}, Geom. Top. 20(4), 1909--1939, (2016).

\bibitem[DL]{Damanik-Lenz} D. Damanik, D. Lenz \emph{Linear repetitivity, I. uniform subadditive ergodic theorems and applications}, Disc. Comp. Geom. 26, 411--428, (2001).

\bibitem[Da]{Danzer} L. Danzer, \emph{Inflation species of planar tilings which are not of locally finite complexity}, Proc. Steklov. Inst. Math. 230, 118--126, (2002).






\bibitem[Du]{Du} F. Durand, \emph{Linearly recurrent subshifts have a finite number of
non-periodic subshift factors}, Ergo. Theo. Dynam. Sys. 20(4), 1061--1078, (2000).


\bibitem[Fr]{Frank} N. P. Frank, \emph{Tilings with infinite local complexity}, in \underline{Mathematics of aperiodic order}, eds. J. Kellendonk, D. Lenz, J. Savinien, Progr. Math., 309, Birkh\"auser/Springer, Basel, 223--257, (2015).

\bibitem[FrRo]{Frank-Robinson} N. P. Frank, E. A. Robinson, Jr., \emph{Generalized $\beta$-expansions, substitution tilings, and local finiteness}, Trans, Amer. Math. Soc. 360, 1163--1177, (2008).

\bibitem[FrS1]{Frank-Sadun1} N. P. Frank, L. Sadun \emph{Topology of some tiling spaces without finite local complexity}, Disc. Cont. Dynam. sys. 23(3), 847--865, (2009).

\bibitem[FrS2]{Frank-Sadun2} N. P. Frank, L. Sadun \emph{Fusion tilings without finite local complexity}, Top. Proc. 43, 235--276, (2014). 


\bibitem[FrRi]{Frettloh-Richard} D. Frettl\"oh, C. Richard, \emph{Dynamical properties of almost repetitive Delone sets}, Disc. Cont. Dynam. Sys. 34(2), 531--556, (2014). 



\bibitem[Ga]{Garber} A. Garber, \emph{On equivalence classes of separated nets}, Modelirovanie i Analiz Informatsionnykh Sistem 16(2), 109--118, (2009).

\bibitem[Gr]{Gr} M. Gromov, \underline{Asymptotic invariants of infinite groups, in Geometric Group Theory}, Vol. 2 (Sussex, 1991), Cambridge Univ. Press, Cambridge,
(1993).



\bibitem[HKeWe]{HKW} A. Haynes, M. Kelly, B. Weiss, \emph{Equivalence relations on separated nets arising from linear toral flows}, Proc. Lond. Math. Soc. 109(5), 1203--1228, (2014).


\bibitem[HKoWa]{Haynes-Koivusalo-Walton} A. Haynes, H. Koivusalo, J. Walton, \emph{A characterization of linearly repetitive cut-and-project sets}, Nonlinearity 31, 515--539, (2018).



 
\bibitem[KW]{Koivusalo-Walton} H. Koivusalo, J. Walton, \emph{Cut and project sets with polytopal window II: linear repetitivity}, Trans. Amer. Math. Soc. 375(7), 5097--5149, (2022).
 
\bibitem[L]{Laczk} M. Laczkovich, \emph{Uniformly spread discrete sets in $\R^d$}, J. Lond. Math. Soc. 46(2), 39--57, (1992).

\bibitem[LP]{Lagarias-Pleasants} J. C. Lagarias, P. A. B. Pleasants, \emph{Repetitive Delone sets and quasicrystals}, Ergo. Theo. Dynam. Sys. 23, 831--867, (2003).

\bibitem[LS]{Lee-Solomyak} J. Y. Lee, B. Solomyak, \emph{On substitution tilings and Delone sets without finite local compexity}, Disc. Cont. Dynam. Sys. 39(6), 3149--3177, (2019). 


\bibitem[Mag]{Magazinov} A. N. Magazinov, \emph{The family of bi-Lipschitz classes of Delone sets in Euclidean space has the cardinality of the continuum}, Proc. of the Steklov. Inst. of Math. 275, 87--98, (2011).

\bibitem[Mat]{Mattila} P. Mattila, \underline{Geometry of sets and measures in Euclidean spaces. Fractals and rectifiability}, volume 44 of Cambridge Studies in Advanced Mathematics. Cambridge University Press, Cambridge, (1995).

\bibitem[McM]{McMullen} C. T. McMullen, \emph{Lipschitz maps and nets in Euclidean space}, Geom. Func. Anal. 8(2), 304--314, (1998).

\bibitem[N]{Navas} A. Navas, \emph{Une remarque \'a propos de l'\'equivalence bilipschitzienne entre des ensembles de Delone (A remark concerning bi-Lipschitz equivalence of Delone sets)}, Comptes Rendus Mathematique, 354(10), 976--979,  (2016).




\bibitem[Sa]{Generalized_Pinwheel} L. Sadun, \emph{Some generalizations of the Pinwheel tiling}, Disc. Comp. Geom. 20(1), 79--110, (1998). 


\bibitem[Sm1]{Yotam Kakutani} Y. Smilansky,
{\em Uniform distribution of Kakutani partitions generated by substitution schemes}, Israel J. Math. 240, 667--710, (2020). 

\bibitem[Sm2]{Frequencies} Y. Smilansky,
{\em Statistics and gap distributions in random Kakutani partitions and multiscale substitution tilings}, J. Math. Anal. Appl. 516(2), 126535, (2022).

\bibitem[SS1]{Smilansky-Solomon1} Y. Smilansky, Y. Solomon, \emph{Multiscale substitution tilings}, Proc. Lond. Math. Soc. 123(6), 517--564, (2021). 

\bibitem[SS2]{Smilansky-Solomon2} Y. Smilansky, Y. Solomon, \emph{A dichotomy for bounded displacement equivalence of Delone sets}, Ergo. Theo. Dynam. Sys. 42(8), 2693--2710, (2022).


\bibitem[yS1]{Solomon11} Y. Solomon, \emph{Substitution tilings and separated nets with similarities to the integer lattice}, Israel J. Math. 181, 445--460, (2011). 

\bibitem[yS2]{Solomon14} Y. Solomon, \emph{A simple condition for bounded displacement}, J. Math. Anal. Appl. 414(1), 134--148, (2014).



\bibitem[bS]{Solomyak-recognizability} B. Solomyak, \emph{Nonperiodicity implies unique composition for self-similar translationally finite tilings}, Disc. Comp. Geom. 20, 265--279, (1998).

\end{thebibliography}
\end{document}